\documentclass[10pt]{amsart}
\usepackage{mathptmx}
\usepackage{amsmath}
\usepackage{amssymb}
\usepackage{array}
\usepackage{geometry}
\usepackage[bookmarks=true,colorlinks=true, pdfstartview=FitV, linkcolor=black, citecolor=blue, urlcolor=black]{hyperref}

\usepackage{color}
\definecolor{DarkRed}{rgb}{0.55,.00,0.2}
\definecolor{DarkGrey}{rgb}{0.35,.35,0.35}

\theoremstyle{definition}

\theoremstyle{remark}

\numberwithin{equation}{section}

\begin{document}

\title{ Index transforms with the square of Bessel functions}

\author{S. Yakubovich}
\address{Department of Mathematics, Faculty of Sciences,  University of Porto,  Campo Alegre str.,  687; 4169-007 Porto,  Portugal}
\email{ syakubov@fc.up.pt}


\date{\today}
\maketitle

\markboth{\rm \centerline{ S.  Yakubovich}}{}
\markright{\rm \centerline{Index Transforms with the square of Bessel functions}}

\begin{abstract}   New index transforms, involving the square of  Bessel  functions  of the first  kind as the kernel  are  considered.   Mapping properties such as the boundedness and  invertibility are investigated for these  operators  in the Lebesgue  spaces.  Inversion theorems are proved.    As an interesting application,  a solution to the initial   value problem for the third  order partial differential equation, involving the Laplacian,  is obtained.  
\end{abstract}

\section{Introduction and preliminary results}

Let $f(x),  \ g(\tau),\  x \in \mathbb{R}_+,\  \tau \in  \mathbb{R}$ be  complex -valued functions.  The main goal of this paper is to investigate mapping properties of the following index transforms   \cite{yak}, involving  the square of Bessel's  function of the first   kind  in  the kernel,  namely,
$$(Ff) (\tau) = {\sqrt\pi\over  \cosh(\pi\tau)} \int_0^\infty  {\rm Re} \left[ J^2_{i\tau} \left( \sqrt x \right) \right]  f(x)dx,   \quad    \tau \in \mathbb{R}, \eqno(1.1)$$
$$(G g) (x) =  \sqrt\pi \   \int_{-\infty}^\infty   {\rm Re}  \left[ J^2_{i\tau} \left(\sqrt  x \right)\right] \   { g(\tau) d\tau\over \cosh(\pi\tau)} ,   \quad   x  \in \mathbb{R}_+, \eqno(1.2)$$
where $i$ is the imaginary unit and ${\rm Re}$ denotes the real part of a complex -valued function.  Bessel's  function of the first kind  $J_\nu(z)$ \cite{erd}, Vol. II  satisfies the differential equation
$$  z^2{d^2u\over dz^2}  + z{du\over dz} + (z^2- \nu^2)u = 0.\eqno(1.3)$$
It has the asymptotic behavior 
$$ J_\nu(z) = \sqrt{2\over \pi z} \cos \left( z- {\pi\over 4} (2\nu+1)\right)  [1+ O(1/z)], \ z \to \infty,\  - {\pi\over 2} < \arg z < {3\pi\over 2}, \eqno(1.4)$$
$$J_\nu(z) = O( z^{\nu} ), \ z \to 0.\eqno(1.5)$$
and the  following series expansion  
$$J_{\nu}(z)= \sum_{k=0}^\infty {(-1)^k (z/2)^{2k+\nu} \over k! \Gamma(k+\nu+1)},\   z,  \nu \in \mathbb{C},\eqno(1.6)$$
where $\Gamma(z)$ is Euler's gamma function \cite{erd}, Vol. I.    Using  the reduction  formula for the gamma function \cite{erd}, Vol. I we find for ${\rm Re} \nu \ge 0$
$$| \Gamma(k+\nu+1) |= | \Gamma(\nu+1) (1+\nu)(2+\nu)\dots (k+\nu)| \ge k!  | \Gamma(\nu+1)|.$$
Hence  we derive from (1.6)
$$|J_{\nu}(z)| \le   e^{-  {\rm Im} \nu \arg z } \sum_{k=0}^\infty {(|z|/2)^{2k+  {\rm Re} \nu} \over k! |\Gamma(k+\nu+1)|} \le 
e^{-  {\rm Im} \nu \arg z } \  { \left(|z| /2\right)^{  {\rm Re} \nu} \over  | \Gamma(\nu+1)| } \sum_{k=0}^\infty {(|z|/2)^{2k} \over (k!)^2}$$
$$\le \  e^{|z| -  {\rm Im} \nu \arg z } \  { \left(|z| /2\right)^{  {\rm Re} \nu} \over  | \Gamma(\nu+1)| } ,$$
coming up to the following inequality for the  Bessel function of the first kind 
$$|J_{\nu}(z)| \le \  \left({|z|\over 2}\right)^{  {\rm Re} \nu} \  \frac {e^{|z| -  {\rm Im} \nu \arg z }}{ | \Gamma(\nu+1)| },\ z,  \nu \in \mathbb{C} .\eqno(1.7)$$
Taking into account the value
$$|\Gamma(1+i\tau)|=  \sqrt {{\pi\tau\over \sinh(\pi\tau)}},$$
inequality (1.7) takes the form
$$ |J_{i\tau}(x)| \le \  e^x\  \sqrt{{\sinh(\pi\tau)\over \pi\tau}},\   x >0,  \tau  \in \mathbb{R}.\eqno(1.8)$$
In the meantime,  appealing to relation (8.4.19.17) in \cite{prud}, Vol. III, we find the following Mellin-Barnes integral representation for the kernel 
$$\Phi_\tau(x)= {\sqrt\pi\over  \cosh(\pi\tau)} \   {\rm Re} \left[ J^2_{i\tau} \left( \sqrt x \right) \right]$$
in (1.1), (1.2),  namely, 
$$ \Phi_\tau(x) = {1\over 2\pi i} \int_{\gamma-i\infty}^{\gamma +i\infty} \frac {\Gamma(s+ i\tau)\Gamma(s-i\tau) \Gamma(1/2-s)}{\Gamma(s) \Gamma (1-s) \Gamma(1- s) } x^{-s} ds, \ x >0,\eqno(1.9)$$
where $\gamma$ is taken from the interval $(0, 1/4)$. It  guarantees the absolute convergence of the integral (1.9)   
since appealing to the Stirling asymptotic formula for the gamma- function \cite{erd}, Vol. I
$$\frac {\Gamma(s+ i\tau)\Gamma(s-i\tau) \Gamma(1/2-s)}{\Gamma(s) \Gamma (1-s) \Gamma(1- s) }  = 
O \left(|s|^ {2\gamma - 3/2}\right),\  |s| \to \infty.\eqno(1.10)$$
Integral  (1.9) can be involved  to represent  $ \Phi_\tau(x)$  for all $x >0,\ \tau \in \mathbb{R}$ as the Fourier cosine transform \cite{tit}.   Indeed, we have

{\bf Lemma 1}. {\it Let $x >0,  \tau \in \mathbb{R}$. Then $ \Phi_\tau(x)$ has the representation in terms of the Fourier cosine transform 
$$ \Phi_\tau(x)  =  -\  {1 \over \sqrt{\pi}  } {\partial \over \partial x}  \int_0^\infty   {\sqrt x \ \cos(\tau u)\over  \cosh (u/2)}  \   {\bf L}_1 \left( 2i \sqrt {x} \  \cosh \left({u\over 2}\right) \right) du,\eqno(1.11) $$
where ${\bf L}_1 (z)$ is the modified Struve function of the index one  \cite{erd}, Vol. II.}

\begin{proof}      Integrating in (1.9) with respect to $x$,  we change the order of integration in the right-hand side of the obtained equality by virtue of the absolute convergence of the iterated integral. As a result we obtain
$$ {1\over x} \int_0^x \Phi_\tau(y) dy  = {1\over 2\pi i} \int_{\gamma-i\infty}^{\gamma +i\infty} \frac {\Gamma(s+ i\tau)\Gamma(s-i\tau) \Gamma(1/2-s)}{\Gamma(s) \Gamma (1-s) \Gamma(2- s) } x^{-s} ds, \ x >0,\eqno(1.12)$$
Hence,  appealing  to the reciprocal formulae via the Fourier cosine transform (cf. formula (1.104) in \cite{yak}) 
$$\int_0^\infty  \Gamma\left(s + i\tau\right)  \Gamma\left(s - i\tau \right)  \cos( \tau y) d\tau
= {\pi\over 2^{2s}}  {\Gamma(2s) \over \cosh^{2s}(y/2)},\ {\rm Re}\ s > 0,\eqno(1.13)$$
$$  \Gamma\left(s + i\tau \right)  \Gamma\left(s - i\tau\right)  
=   { \Gamma(2s)  \over 2^{2s-1}}  \int_0^\infty   {\cos(\tau y)  \over \cosh^{2s} (y/2)} \ dy,\eqno(1.14)$$ 
we replace  the gamma-product $ \Gamma\left(s + i\tau \right)  \Gamma\left(s - i\tau\right)$ in the integral (1.12)  by its integral representation (1.14) and change the order of integration via the absolute convergence, which can be justified using the Stirling asymptotic formula for the gamma-function.  Then, employing the duplication formula for the gamma- function \cite{erd}, Vol. I,   we derive 
$$ {1\over x} \int_0^x \Phi_\tau(y) dy  = {1\over 2 \pi\sqrt \pi  i} \int_0^\infty   \cos(\tau u)  \int_{\gamma-i\infty}^{\gamma +i\infty} \frac {\Gamma(s+ 1/2)  \Gamma(1/2-s)}{ \Gamma (1-s) \Gamma(2- s) } \left( x \cosh^2(u/2) \right)^{-s}  ds du  $$
$$= -\  {1 \over \sqrt{\pi x}  }  \int_0^\infty   {\cos(\tau u)\over \cosh (u/2)}  \   {\bf L}_1 \left( 2i \sqrt {x} \   \cosh \left({u\over 2}\right) \right) dy,$$
where the  inner integral with respect to $s$ is calculated via Slater's theorem and relation (7.14.2.79) in \cite{prud}, Vol. III.    Hence after the differentiation with respect to $x$  we arrive at (1.11),  completing the proof of Lemma 1. 

\end{proof} 

Further,  recalling  the Mellin-Barnes representation (1.9) of the kernel  $\Phi_\tau(x)$, we will derive an ordinary differential equation whose particular solution is  $\Phi_\tau(x)$. Precisely,   it is given by 

{\bf Lemma 2}.    {\it For each $\tau \in \mathbb{R}$ the function  
$$\Phi_\tau(x)= {\sqrt\pi\over  \cosh(\pi\tau)} \   {\rm Re} \left[ J^2_{i\tau} \left( \sqrt x \right) \right]$$
is a fundamental solution of the following third order differential equation with variable  coefficients}
$$x^2 {d^3 \Phi_{\tau} \over dx^3} + 3x  {d^2 \Phi_\tau \over dx^2}  +  \left( \tau^2+ x- 7 \right)  {d\Phi_\tau \over dx}  + {\Phi_\tau \over 2} = 0,\ x >0. \eqno(1.15)$$

\begin{proof}  The asymptotic behavior (1.10) of the integrand in (1.9) and the absolute and uniform convergence of the integral allow to integrate (1.9) twice  with respect to $x$.  Hence making use the reduction formula for the gamma- function \cite{erd}, Vol. I   we obtain 
$$\int_0^x du \int_0^u  \Phi_{\tau}( v) dv =  {1\over 2\pi i} \int_{\gamma-i\infty}^{\gamma +i\infty} \frac {\Gamma(s+ i\tau)\Gamma(s-i\tau) \Gamma(1/2-s)}{\Gamma(s) \Gamma (1-s) \Gamma(3- s) } x^{ 2 -s} ds.\eqno(1.16)$$
Then, plainly,  with the reduction formula for the gamma-function and simple changes of variables we get  
$$\left( x {d\over dx} \right)^2  \left[ {1\over x^2}  \int_0^x du \int_0^u  \Phi_{\tau}( v) dv \right] 
= {1\over 2\pi i} \int_{\gamma-i\infty}^{\gamma +i\infty} \frac { s^2 \Gamma(s+ i\tau)\Gamma(s-i\tau) \Gamma(1/2-s)}{\Gamma(s) \Gamma (1-s) \Gamma(3- s) } x^{  -s} ds $$

$$= -{\tau^2\over x^2 }  \int_0^x du \int_0^u  \Phi_{\tau}( v) dv +  {1\over 2\pi  i} \int_{\gamma-i\infty}^{\gamma +i\infty} \frac {\Gamma(1+ s+ i\tau)\Gamma(1+ s-i\tau) \Gamma(1/2-s)}{\Gamma(s)\Gamma (1-s) \Gamma (3-s)} x^{-s} ds $$

$$= -{\tau^2\over x^2 }  \int_0^x du \int_0^u  \Phi_{\tau}( v) dv -   {1\over 2\pi  i} \int_{1+ \gamma-i\infty}^{1+ \gamma +i\infty} \frac {\Gamma( s+ i\tau)\Gamma( s-i\tau) \Gamma(1/2-s) (1/2-s) }{\Gamma(s)\Gamma (1-s) \Gamma (1-s)(3-s)(2-s) (1-s) } x^{1 -s} ds. $$

Hence after multiplication by $x^2$, we differentiate three times the obtained equality  and use (1.9) to find   

$$ {d^3\over dx^3}  \left [x^2 \left( x {d\over dx} \right)^2  \left[ {1\over x^2}  \int_0^x du \int_0^u  \Phi_{\tau}( v) dv \right] \right]  =   - \tau^2 {d\Phi_\tau \over dx}  - \sqrt x {d\over dx}  \left( \sqrt x \ \Phi_\tau (x)\right)$$
$$=   - \left( \tau^2+ x\right)  {d\Phi_\tau \over dx}  - {\Phi_\tau \over 2}.\eqno(1.17) $$
Meanwhile,   the left-hand side of (1.17) can be simplified  as follows  
$$ {d^3\over dx^3}  \left [x^2 \left( x {d\over dx} \right)^2  \left[ {1\over x^2}  \int_0^x du \int_0^u  \Phi_{\tau}( v) dv \right] \right] = {d^3\over dx^3}  \left [x^3 {d\over dx}   \left[ - {2\over x^2}  \int_0^x du \int_0^u  \Phi_{\tau}( v) dv  + {1\over x} \int_0^x  \Phi_{\tau}( v) dv \right] \right] $$

$$=  {d^3\over dx^3}  \left [ - 4  \int_0^x du \int_0^u  \Phi_{\tau}( v) dv  - 3x    \int_0^x  \Phi_{\tau}( v) dv+  x^2 \Phi_{\tau} \right]  =  x^2 {d^3 \Phi_\tau \over dx^3}  +3x  {d^2 \Phi_\tau \over dx^2}  - 7  {d \Phi_\tau \over dx} .$$

Finally, combining with (1.17), we arrive at the equation (1.15), completing the proof of Lemma 2.
\end{proof}

\section {Boundedness  and inversion properties of the index transform (1.1)}

In order to investigate the mapping properties of the index transform (1.1) we will use the Mellin transform technique developed in \cite{yal}.   Precisely, the Mellin transform is defined, for instance, in  $L_{\nu, p}(\mathbb{R}_+),\ 1 < p \le 2$ (see details in \cite{tit}) by the integral  
$$f^*(s)= \int_0^\infty f(x) x^{s-1} dx,\eqno(2.1)$$
 being convergent  in mean with respect to the norm in $L_q(\nu- i\infty, \nu + i\infty),\   q=p/(p-1)$.   Moreover, the  Parseval equality holds for $f \in L_{\nu, p}(\mathbb{R}_+),\  g \in L_{1-\nu, q}(\mathbb{R}_+)$
$$\int_0^\infty f(x) g(x) dx= {1\over 2\pi i} \int_{\nu- i\infty}^{\nu+i\infty} f^*(s) g^*(1-s) ds.\eqno(2.2)$$
The inverse Mellin transform is given accordingly
 $$f(x)= {1\over 2\pi i}  \int_{\nu- i\infty}^{\nu+i\infty} f^*(s)  x^{-s} ds,\eqno(2.3)$$
where the integral converges in mean with respect to the norm  in   $L_{\nu, p}(\mathbb{R}_+)$
$$||f||_{\nu,p} = \left( \int_0^\infty  |f(x)|^p x^{\nu p-1} dx\right)^{1/p}.\eqno(2.4)$$
In particular, letting $\nu= 1/p$ we get the usual space $L_1(\mathbb{R}_+)$.  Further, denoting by $C (\mathbb{R}_+)$ the space of bounded continuous functions, we have 

{\bf Theorem 1.}   {\it The index transform  $(1.1)$  is well-defined as a  bounded operator $F: L_1\left(\mathbb{R}_+; e^{2\sqrt x} dx\right) \to C (\mathbb{R})$ and the following norm inequality takes place
$$||Ff||_{C (\mathbb{R})} \equiv \sup_{\tau \in \mathbb{R}} | (Ff)(\tau)| \le \sqrt \pi ||f||_{L_1\left(\mathbb{R}_+; \ e^{2\sqrt x} dx\right)}.\eqno(2.5)$$
Moreover,  let  the Mellin transform $(2.1)$ of $f$ satisfy  the condition
$${f^*(s)\over \Gamma (s)\Gamma(1-s)} \in L_q (1-\nu-i\infty, 1-\nu+ i\infty),\eqno(2.6)$$
where
$$\  0 < \nu < {3\over 4} - {1\over 2p},\   1 < p < 2,\   {1\over p} + {1\over q} =1.\eqno(2.7)$$
Then 
$$(Ff)(\tau)=  {\sqrt\pi\over  \cosh(\pi\tau)} \int_0^\infty  {\rm Re} \left[ I_{i\tau} \left({x\over 2} \right) \right]  e^{-x/2}  \varphi (x) dx,\eqno(2.8)$$
where $I_\mu(z)$ is the modified Bessel function of the first kind, 
$$\varphi(x)= {1\over 2\pi i}  \int_{1-\nu- i\infty}^{1-\nu+i\infty} {f^*(s)\over \Gamma(s)}   x^{-s} ds,\eqno(2.9)$$
and  integral $(2.9)$ converges with respect to the norm in $L_{1-\nu, p}(\mathbb{R}_+)$ .}

\begin{proof} The proof of the norm inequality (2.5) is straightforward from  (1.1) and  inequality (1.8).  The continuity of $(Ff ) (\tau)$ follows from the absolute and uniform convergence of the corresponding integral. In fact, we derive
$$ | (Ff)(\tau)| \le {\sqrt\pi \over  \cosh(\pi\tau)} \int_0^\infty \left| J _{i\tau} \left(\sqrt x \right) \right|^2   |  f(x)| dx $$
$$\le \sqrt \pi \   {\tanh (\pi\tau) \over  \pi\tau } \int_0^\infty     |  f(x)|\  e^{2\sqrt x} \  dx 
 \le \sqrt \pi ||f||_{L_1\left(\mathbb{R}_+; \ e^{2\sqrt x} dx\right)}.$$
This proves (2.5).   Next,   $f^*(1-s)  \in L_q(\nu-i\infty,  \nu+ i\infty)$ via condition (2.6).   Indeed, we have
$$\int_{\nu-i\infty}^{ \nu  +i\infty} |f^*(1-s) |^q  |ds| =  \int_{1-\nu-i\infty}^{1- \nu  +i\infty} |\Gamma(s) \Gamma (1-s) f^*(s) |^q {|ds|\over |\Gamma(s) \Gamma(1-s)|^q} $$ 
$$ \le \left[\Gamma(\nu)\Gamma(1-\nu) \right]^q    \int_{1-\nu-i\infty}^{1- \nu  +i\infty} |f^*(s) |^q {|ds|\over |\Gamma(s) \Gamma(1-s)|^q} < \infty.\eqno(2.10)$$ 
Therefore,  via Theorem 86 in \cite{tit} we  see that $f \in L_{1-\nu,p}(\mathbb{R}_+)$.  Meanwhile, the asymptotic behavior (1.10)  of  the  integrand in (1.9) guarantees that it belongs to the space $L_p (\nu-i\infty,  \nu+ i\infty) $ with $\nu$, satisfying condition (2.7).  Consequently,  (1.9) and the Parseval identity (2.2) lead to the following representation 
$$(Ff)(\tau)= {1\over 2\pi  i} \int_{\nu-i\infty}^{\nu  +i\infty} \frac {\Gamma(s+ i\tau)\Gamma(s-i\tau) \Gamma(1/2-s)}{\Gamma(s)\Gamma^2 (1-s)} f^*(1-s) x^{-s} ds.\eqno(2.11)$$
Analogously to (2.10) we show that  $f^*(1-s)  /\Gamma(1-s) \in L_q(\nu-i\infty,  \nu+ i\infty)$ and hence (see (2.9)) $\varphi \in  L_{1-\nu, p}(\mathbb{R}_+).$   But from the asymptotic behavior of the modified Bessel function of the first kind (see \cite{erd}, Vol. II) we find 
$$e^{-x/2} {\rm Re} \left[ I_{i\tau} \left({x\over 2}\right) \right] = O(1),\ x \to 0,$$
$$e^{-x/2} {\rm Re} \left[ I_{i\tau} \left({x\over 2}\right) \right] = O\left( {1\over \sqrt x}\right),\ x \to \infty.$$
Therefore,   the kernel of (2.8) belongs to $L_{\nu, q}(\mathbb{R}_+)$  and integral  converges absolutely.  Moreover, Lemma 1 in \cite{skal},  the Parseval identity  (2.2) and (2.11) establish equality (2.8), completing  the proof of Theorem 1.   

\end{proof}

The inversion formula for the transform (1.1) is  given by 

{\bf Theorem 2.}  {\it Under conditions of Theorem 1 let also the Mellin transform $f^*(s)$ be analytic in the strip 
$${1\over 2p}- {1\over 4}  < {\rm Re}\  s = 1-\nu  <  {7\over 4} - {1\over 2p},\ 1 < p < 2 $$ 
and 
$${f^*(s)\over \Gamma(s) \Gamma (1-s)} \in L_q (1-\nu-i\infty, 1-\nu+ i\infty) \cap L_1 (1-\nu-i\infty, 1-\nu+ i\infty),\ |\nu| < {3\over 4} - {1\over 2p}. \eqno(2.12)$$
Then, assuming that the index transform $(1.1)$ satisfies the integrability condition $(Ff)(\tau) \in L_1(\mathbb{R}_+; \tau e^{\pi\tau} d\tau)$, it has the following inversion formula for all $x >0$ 

$$f(x)=  2 \sqrt\pi {d\over dx}   \int_0^\infty \tau  \coth(\pi\tau)\left[ {  {\rm Im}  \left[ J^2_{i\tau} (\sqrt x)  \right] \over \sinh (\pi\tau)}  -  {1\over \pi}  \  \left[ J_{i\tau} (\sqrt x)  {\partial  J_{\varepsilon - i\tau} (\sqrt x) \over \partial \varepsilon }  \biggr\rvert_{\varepsilon=0} \right.\right.$$
$$\left.\left.   +  J_{ -i\tau} (\sqrt x)      {\partial  J_{\varepsilon +i\tau} (\sqrt x) \over \partial \varepsilon }\biggr\rvert_{\varepsilon=0}  \right] \right]   (Ff)(\tau) d\tau,\eqno(2.13)$$
where ${\rm Im} $ denotes the imaginary part of a complex -valued function and the corresponding  integral converges absolutely.}

\begin{proof}   Following the same scheme as in \cite{skal}, Th. 2, we verify conditions of the Lebedev expansion theorem for the index transform with the square of the Macdonald function as the kernel  $K^2_{i\tau} (\sqrt x)$ \cite{erd}, Vol. II, coming up to the equality 
$$\int_x^\infty    {1\over 2\pi i}  \int_{1-\nu- i\infty}^{1-\nu+i\infty} {f^*(s)\   y^{-s} \over \Gamma(s) \Gamma(1-s)}   ds  dy = {2\over \pi^2\sqrt\pi } \int_0^\infty \tau \sinh(2\pi\tau) K^2_{i\tau} (\sqrt x) (Ff)(\tau) d\tau,\ x >0.\eqno(2.14)$$
Fulfilling the integration in the left-hand side of (2.14) when $\nu \in (-3/4 + 1/(2p),\ 0)$,  it becomes 
$$-  {1\over 2\pi i}  \int_{1-\nu- i\infty}^{1-\nu+i\infty} {f^*(s)\   x^{ 1 -s} \over \Gamma(s) \Gamma(2-s)}   ds 
= {2\over \pi^2\sqrt\pi } \int_0^\infty \tau \sinh(2\pi\tau) K^2_{i\tau} (\sqrt x) (Ff)(\tau) d\tau.\eqno(2.15)$$
Hence
$$-  {1\over 2\pi i}  \int_0^\infty {dy\over (x+y)^2 }  \int_{1-\nu- i\infty}^{1-\nu+i\infty} {f^*(s)\   y^{ 1 -s} \over \Gamma(s) \Gamma(2-s)}   ds $$
$$=   {2\over \pi^2\sqrt\pi }  \int_0^\infty {dy\over (x+y)^2 }  \int_0^\infty \tau \sinh(2\pi\tau) K^2_{i\tau} (\sqrt y) (Ff)(\tau) d\tau.$$
The interchange of the order of integration in both sides of the latter equality is allowed due to Fubini's theorem under imposed conditions.  Thus, taking into account the value of the elementary beta-integral and the inverse Mellin transform (2.3), we end up with the inversion formula of the index transform (1.1), namely
$$ f(x)=  {2\over \pi^2\sqrt\pi }   \int_0^\infty \tau \sinh(2\pi\tau) K(x,\tau) (Ff)(\tau) d\tau, \ x >0, \eqno(2.16)$$
where
$$ K(x,\tau) =  -  \int_0^\infty { K^2_{i\tau} (\sqrt y) \ dy\over (x+y)^2 }.\eqno(2.17)$$  
In the meantime  integral  (2.17), which is absent in the literature,  can be calculated in the following way.   
Using relations (8.4.23.27) and (8.4.2.5) in \cite{prud}, Vol. III and the Parseval equality (2.2), we write
$$ -  \int_0^\infty { K^2_{i\tau} (\sqrt y) \ dy\over (x+y)^2 } =  {\sqrt\pi \over 4\pi   i}  {d\over dx} \int_{\gamma -i\infty}^{\gamma   +i\infty} \Gamma(s+ i\tau)\Gamma(s-i\tau) \Gamma^2(s) \Gamma (1-s)  {x^{-s}\over \Gamma(1/2+ s) }ds,\eqno(2.18)$$
where $0 < \gamma < 1$.   On the other hand, the Lebedev inequality for the Macdonald function (cf. \cite{yal}, p.  99)
$$|K_{i\tau} (x)| \le {x^{-1/4}\over \sqrt {\sinh(\pi\tau)}},\ x, \tau > 0\eqno(2.19)$$
and the estimate 
$$\int_0^\infty \tau \sinh(2\pi\tau) \left| K(x,\tau) (Ff)(\tau)\right|  d\tau \le  2 \int_0^\infty \tau \cosh(2\pi\tau) \left|(Ff)(\tau)\right|  d\tau \int_0^\infty  {  y^{-1/4}\ dy\over (x+y)^2 } < \infty, x >0$$
under condition $(Ff)(\tau)  \in L_1 (\mathbb{R}_+; \tau e^{\pi\tau} d\tau )$ permit us to write (2.16) in the form
$$f(x)=   {2\over \pi^2\sqrt\pi }  {d\over dx}  \int_0^\infty \int_0^\infty \tau \sinh(2\pi\tau)  K^2_{i\tau} (\sqrt y)  (Ff)(\tau) { dy \ d\tau\over x+y }, \eqno(2.20)$$
where via the uniform convergence and (2.18) 
$$ \int_0^\infty { K^2_{i\tau} (\sqrt y) dy\over x+y } = {\sqrt\pi \over 4\pi   i}  \int_{\gamma -i\infty}^{\gamma   +i\infty} \Gamma(s+ i\tau)\Gamma(s-i\tau) \Gamma^2(s) \Gamma (1-s)  {x^{-s}\over \Gamma(1/2+ s) }ds $$
$$=  {\sqrt\pi \over 4\pi   i}  \lim_{\varepsilon \to 0} \int_{\gamma -i\infty}^{\gamma   +i\infty} \Gamma(s+ i\tau)\Gamma(s-i\tau) \Gamma(s+\varepsilon)  \Gamma(s- \varepsilon) \Gamma (1-s)  {x^{-s}\over \Gamma(1/2+ s) }ds. \eqno(2.21)$$
Hence the Slater theorem is applied, involving simple left-hand poles $s= \pm i\tau -n,\ s= \pm \varepsilon - n,\  n \in \mathbb{N}_0$.  Calculating the corresponding residues and then  passing to the limit when $\varepsilon \to 0$  with the use of  relation (7.15.1.3) in \cite{prud}, Vol. III, we find the result 
$$ {1\over 2\pi   i}   \lim_{\varepsilon \to 0}  \int_{\gamma -i\infty}^{\gamma   +i\infty} \Gamma(s+ i\tau)\Gamma(s-i\tau) \Gamma(s+\varepsilon)  \Gamma(s- \varepsilon) \Gamma (1-s)  {x^{-s}\over \Gamma(1/2+ s) }ds$$
$$=   \pi^2\sqrt\pi \left[ {  {\rm Im}  \left[ J^2_{i\tau} (\sqrt x)  \right] \over \sinh^3(\pi\tau)}  +   \lim_{\varepsilon \to 0} 
\  \frac{J_{-\varepsilon -i\tau} (\sqrt x)    J_{-\varepsilon +i\tau} (\sqrt x)- J_{\varepsilon -i\tau} (\sqrt x)    J_{\varepsilon +i\tau} (\sqrt x)}  {\sin(\pi\varepsilon) \left( \cosh(2\pi\tau)- \cos(2\pi\varepsilon)\right)  } \right] $$
$$ =   {\pi^2\sqrt\pi\over \sinh^2(\pi\tau)} \left[ {  {\rm Im}  \left[ J^2_{i\tau} (\sqrt x)  \right] \over \sinh (\pi\tau)}  -  {1\over \pi}  
\  \left[ J_{i\tau} (\sqrt x)  {\partial  J_{\varepsilon - i\tau} (\sqrt x)\over \partial \varepsilon }  \biggr\rvert_{\varepsilon=0}  +  J_{ -i\tau} (\sqrt x)      {\partial J_{\varepsilon +i\tau} (\sqrt x) \over \partial \varepsilon } \biggr\rvert_{\varepsilon=0}  \right] \right] .$$
Therefore, combining with (2.20), (2.21),  we arrive at the inversion formula (2.13),    completing  the proof of Theorem 2. 

\end{proof}

\section{The index transform (1.2)} 

The boundedness and invertibility properties for the index transform  (1.2) will be examined below.    We begin with 

{\bf Theorem 3.}  {\it  The $g \in L_1(\mathbb{R})$. Then $(Gg)(x)$ is bounded continuous on $\mathbb{R}_+$ and it holds  
$$\sup_{x >0}   \left|  (Gg)(x) \right|   \le \sqrt \pi ||g ||_{L_1(\mathbb{R})}.\eqno(3.1)$$
Moreover,   if $(Gg) (x) \in L_{\nu,1}(\mathbb{R}_+),\  0 < \nu < 1/2$, then for all $y >0$ }
$${1\over 2\pi i}   \int_{\nu -i\infty}^{\nu  +i\infty} \Gamma(s)\Gamma^2(1-s) (Gg)^*(s)  y^{ -s} ds\
=  \sqrt \pi  \   \int_{-\infty}^\infty  e^{y/2} \   K_{i\tau} \left({y\over 2} \right)  {g(\tau) \over \cosh(\pi\tau) } d\tau. \eqno(3.2)$$

\begin{proof}   First we observe from the familiar Poisson integral for the Bessel function \cite{erd}, Vol. II that for all $x >0$
$$\left| J_{i\tau} (\sqrt x)\right| \le { \sqrt \pi\over {\left| \Gamma (i\tau + 1/2)\right| }}.$$
 Therefore, 
 $$  {\left| {\rm Re}  \left[ J^2_{i\tau} \left(\sqrt  x \right)\right] \right|  \over \cosh(\pi\tau)} \  \le \  1,\eqno(3.3)$$
and  we find  the estimate

$$ \left|  (Gg)(x) \right|  \le \sqrt\pi \   \int_{-\infty}^\infty    |g(\tau)| d\tau =   \sqrt \pi ||g ||_{L_1(\mathbb{R})},$$
which yields  (3.1).   Next, integrating both sides of (1.2) with respect to $x$ and changing the order of integration under  condition $g \in L_1(\mathbb{R})$, we get 
$$\int_0^x (G g) (y)\  dy  =   \int_{-\infty}^\infty   \int_0^x  \Phi_\tau(u) \   g(\tau) du\ d\tau. $$
Hence,   employing (1.12), we derive 
$${1\over x} \int_0^x (G g) (y)\  dy  =  {1\over 2\pi i}  \int_{-\infty}^\infty    g(\tau)  \int_{\nu -i\infty}^{\nu  +i\infty} \frac {\Gamma(s+ i\tau)\Gamma(s-i\tau) \Gamma(1/2-s)}{\Gamma(s) \Gamma (1-s) \Gamma(2- s) } x^{ -s} ds\ d\tau$$
$$= {1\over 2\pi  i}    \int_{\nu -i\infty}^{\nu  +i\infty} \int_{-\infty}^\infty    g(\tau) \frac {\Gamma(s+ i\tau)\Gamma(s-i\tau) \Gamma(1/2-s)}{\Gamma(s) \Gamma (1-s) \Gamma(2- s) } x^{-s}  d\tau\ ds,\eqno(3.4)$$
where the interchange of the order of integration is motivated by the straightforward estimate
$$\int_{\nu -i\infty}^{\nu  +i\infty} \int_{-\infty}^\infty    |g(\tau)|  \left| \frac {\Gamma(s+ i\tau)\Gamma(s-i\tau) \Gamma(1/2-s)}{\Gamma(s) \Gamma (1-s) \Gamma(2- s) } x^{1 -s} \right|  d\tau\ |ds| $$
$$\le  x^{1-\nu } B(\nu, \nu) \int_{-\infty}^\infty    |g(\tau)|  d\tau \  \int_{\nu -i\infty}^{\nu  +i\infty} \left| \frac {\Gamma(2s) \Gamma(1/2-s)}{\Gamma(s) \Gamma (1-s) \Gamma(2- s) } \right|  \  |ds|  < \infty,\  0 < \nu < {1\over 2}$$
and $B(a,b)$ is Euler's beta-function \cite{erd}, Vol. I.     But under the condition  $(Gg) (x) \in L_{\nu,1}(\mathbb{R}_+)$ the Mellin transform (2.1) of the left-hand side in (3.4) exists and we have 
$$\int_0^\infty x^{s-2} \int_0^x (G g) (y)\  dy dx=   \int_0^\infty (Gg)(y) \int_y^\infty x^{s-2} dx dy = {(Gg)^* (s) \over 1-s},\ {\rm Re} s =\nu \in \left(0,\ {1\over 2} \right). $$
Thus (3.4) and the reduction formula for the gamma-function imply
$$ \Gamma(s) \Gamma^2 (1-s) (Gg)^* (s)  =   \Gamma(1/2-s) \int_{-\infty}^\infty    g(\tau) \Gamma(s+ i\tau)\Gamma(s-i\tau)  d\tau.$$
Hence the inverse Mellin transform (2.3) and  relation (8.4.23.5) in \cite{prud}, Vol. III will lead us to (3.2), completing the proof of Theorem 3. 
\end{proof} 

The inversion formula for the index transform (1.2) is given by

{\bf Theorem 4}.  {\it  Let $g(z/i)$ be an even analytic function in the strip $D= \left\{ z \in \mathbb{C}: \ |{\rm Re} z | < \alpha < 1/2\right\} ,\  g(0)=g^\prime (0)=0$ and $g(z/i)$ is absolutely  integrable over any vertical line in  $D$.   If $(Gg) (t) \in L_{1}((0,1);  t^{-1} dt) \cap L_{1}((1, \infty);  dt) $, then for all  $x \in \mathbb{R}$ the  inversion formula holds for the index transform (1.2)} 
$$ g(x)  = \int_0^\infty (Gg)(y) \left[ {\cosh(\pi x ) \over \pi \sqrt\pi}    + \pi \  x y \coth(\pi x)\   {d\over dy}    \left[ {  {\rm Im}  \left[ J^2_{ix } (\sqrt y)  \right] \over \sinh (\pi x )} \right.\right.$$
$$\left.\left.  -  {1\over \sqrt \pi}  \  \left[ J_{i x} (\sqrt y)  {\partial  J_{\varepsilon - i x } (\sqrt y)\over \partial \varepsilon }  \biggr\rvert_{\varepsilon=0}  +  J_{ -i x } (\sqrt y)      {\partial J_{\varepsilon +i x } (\sqrt y) \over \partial \varepsilon } \biggr\rvert_{\varepsilon=0}  \right] \right] \right] {dy\over y}.\eqno(3.5)$$

\begin{proof}    Indeed,  recalling  (3.2), we multiply its both sides by $ e^{-y/2} K_{ix} \left({y/2} \right) y^{\varepsilon -1}$ for some positive $\varepsilon \in (0,1)$ and integrate with respect to $y$ over $(0, \infty)$.  Hence changing the order of integration in the left-hand side of the obtained equality due to the absolute convergence of the  iterated integral, we  appeal  to relation (8.4.23.3) in \cite{prud}, Vol. III to find  
$$ {1\over 2\pi i}   \int_{\nu -i\infty}^{\nu  +i\infty} \frac {\Gamma(\varepsilon- s+ ix)\Gamma(\varepsilon- s-ix)}{ \Gamma(1/2+ \varepsilon -s)}\Gamma(s)\Gamma^2(1-s) (Gg)^*(s) ds$$
$$ =  \ \int_0^\infty  K_{ix} \left({y\over 2} \right) y^{\varepsilon -1} \int_{-\infty}^\infty  K_{i\tau} \left({y\over 2} \right)  {g(\tau) \over \cosh(\pi\tau) } d\tau dy. \eqno(3.6)$$

Meanwhile,  the right-hand side of (3.6) can be treated in the same manner as in the proof of Theorem 4 in \cite{skal}. In fact,   the evenness of $g$, a  representation  of the Macdonald function in terms of the modified Bessel function of the first kind $I_z(y)$ \cite{erd}, Vol. II and a simple substitution give the equality 
$$  \int_0^\infty  K_{ix} \left({y\over 2} \right) y^{\varepsilon -1} \int_{-\infty}^\infty  K_{i\tau} \left({y\over 2} \right)  {g(\tau) \over \cosh(\pi\tau) } d\tau dy$$

$$= 2 \pi i \int_0^\infty  K_{ix} \left({y\over 2} \right) y^{\varepsilon -1} \int_{-i\infty}^{i\infty}   I_{ z} \left({y\over 2} \right)  {g(z/i) \over \sin (2\pi z) } dz\  dy. \eqno(3.7)$$
On the other hand, according to our assumption $g(z/i)$ is analytic in the vertical  strip $0\le  {\rm Re}  z < \alpha< 1/2$,  $g(0)=g^\prime (0)=0$ and integrable  in the strip.  Hence,  appealing to the inequality for the modified Bessel   function of the first  kind  (see \cite{yal}, p. 93)
 $$|I_z(y)| \le I_{  {\rm Re} z} (y) \  e^{\pi |{\rm Im} z|/2},\   0< {\rm Re} z < \alpha,$$
one can move the contour to the right in the latter integral in (3.7). Then 

$$2 \pi i  \int_0^\infty  K_{ix} \left({y\over 2} \right) y^{\varepsilon -1} \int_{-i\infty}^{i\infty}   I_{ z} \left({y\over 2} \right)  {g(z/i) \over \sin (2\pi z) } dz\  dy$$

$$= 2 \pi i    \int_0^\infty  K_{ix} \left({y\over 2} \right) y^{\varepsilon -1} \int_{\alpha -i\infty}^{\alpha + i\infty}   I_{ z} \left({y\over 2} \right)  {g(z/i) \over \sin (2\pi z) } dz\  dy.$$
Now ${\rm Re} z >0$,  and  it is possible to pass to the limit under the integral sign when $\varepsilon \to 0$ and to change the order of integration due to the absolute and uniform convergence.  Therefore the value of the integral (see relation (2.16.28.3) in \cite{prud}, Vol. II)
$$\int_0^\infty K_{ix}(y) I_z(y) {dy\over y} = {1\over x^2 + z^2} $$ 
leads us to the equalities 

$$\lim_{\varepsilon \to 0}  2 \pi i  \int_0^\infty  K_{ix} \left({y\over 2} \right) y^{\varepsilon -1} \int_{-i\infty}^{i\infty}   I_{ z} \left({y\over 2} \right)  {g(z/i) \over \sin (2\pi z) } dz\  dy$$

$$=    2 \pi i  \int_{\alpha -i\infty}^{\alpha + i\infty}   {g(z/i) \over (x^2+ z^2) \sin (2\pi z) } dz =  \pi   i 
 \left( \int_{-\alpha +i\infty}^{- \alpha- i\infty}   +   \int_{\alpha -i\infty}^{ \alpha+  i\infty}   \right)  {  g(z/i) \  dz \over (z-ix) \  z \sin(2\pi z)}. \eqno(3.8)$$
Hence conditions of the theorem allow to apply the Cauchy formula in the right-hand side of the latter equality in (3.8).  Thus 
$$\lim_{\varepsilon \to 0}  2 \pi i   \ \int_0^\infty  K_{ix} \left({y\over 2} \right) y^{\varepsilon -1} \int_{-i\infty}^{i\infty}   I_{ z} \left({y\over 2} \right)  {g(z/i) \over \sin (2\pi z) } dz\  dy =  { 2\pi^{2} \  g(x) \over  x\sinh (2\pi x)} ,\quad x \in \mathbb{R} \backslash \{0\}.\eqno(3.9)$$
Now,  recalling the Parseval identity  (2.2), the left-hand side of (3.6) can be rewritten in the form
$$ {1\over 2\pi i}   \int_{\nu -i\infty}^{\nu  +i\infty} \frac {\Gamma(\varepsilon- s+ ix)\Gamma(\varepsilon- s-ix)}{ \Gamma(1/2+ \varepsilon -s)}\Gamma(s)\Gamma^2(1-s)  (Gg)^*(s) ds = \int_0^\infty (Gg)(y) \Psi_{\varepsilon}  (x,y)  {dy\over y},\eqno(3.10)$$
where
$$\Psi_{\varepsilon}  (x,y) = {1\over 2\pi i}   \int_{1-\nu -i\infty}^{1-\nu  +i\infty} \frac {\Gamma(s+ \varepsilon- 1+ ix)\Gamma(s+ \varepsilon- 1-ix)}{ \Gamma(s+ \varepsilon - 1/2)}\Gamma^2(s) \Gamma(1-s) \  y^{1-s} ds.\eqno(3.11)$$ 
Meanwhile, relations (8.4.23.3),  (8.4.11.3) in \cite{prud}, Vol. III and the Parseval equality (2.2) give the representation (3.11)  accordingly 
$$\Psi_{\varepsilon}  (x,y) =  {y^\varepsilon\over \sqrt\pi} \int_0^\infty \int_1^\infty  \exp \left( u(1-t) - {y\over 2u} \right) K_{ix} \left({y\over 2u} \right)   {dt \  du\over t u^{\varepsilon}}$$
$$=  {y^\varepsilon  \over \sqrt\pi} \int_0^\infty {dt\over 1+t} \int_0^\infty  \exp \left( - {t\over u} \ -\  {yu \over 2}  \right) K_{ix} \left({y u \over 2} \right)   { du\over  u^{2-\varepsilon}} $$
$$=  -  {y^\varepsilon   \over \sqrt\pi} \int_0^\infty {dt\over 1+t} \ {d\over dt} \int_0^\infty  \exp \left( - {t\over u} \ -\  {yu \over 2}  \right) K_{ix} \left({y u \over 2} \right)   u^{\varepsilon -1}  du,\ \varepsilon >0,$$ 
where the interchange of the order of integration and differentiation under the integral sign are permitted owing to the absolute and uniform convergence.  Integrating by parts in the latter integral with respect to $t$ and using relation (2.16.9.3) in \cite{prud}, Vol. II,  we find 
$$\Psi_{\varepsilon}  (x,y) =  {|\Gamma (\varepsilon +ix)|^2\over \Gamma(\varepsilon + 1/2)} -
  {y^\varepsilon   \over \sqrt\pi} \int_0^\infty {dt\over (1+t)^2} \  \int_0^\infty  \exp \left( - {t\over u} \ -\  {yu \over 2}  \right) K_{ix} \left({y u \over 2} \right)   u^{\varepsilon -1}  du.\eqno(3.12)$$
Moreover, it is possible now to pass to the limit under the integral sign in (3.12) when $\varepsilon \to 0$ and then to employ relation (2.16.9.1) in \cite{prud}, Vol. II.  Hence, reminding (2.17), we obtain 

$$\Theta(x,y)= \lim_{\varepsilon \to 0} \Psi_{\varepsilon}  (x,y) =  {\sqrt\pi\over x\sinh(\pi x)}  -
  {2  \over \sqrt\pi} \int_0^\infty { K^2_{ix}\left(\sqrt{yt} \right) dt\over (1+t)^2}= {\sqrt\pi\over x\sinh(\pi x)}  +    2y   K(y,x).\eqno(3.13)$$

In the meantime,  passing to the limit in (3.10) when $\varepsilon \to 0$, we do it under the integral sign in its right-hand side for each $x \in \mathbb{R} \backslash \{0\}$, appealing to the dominated convergence theorem   due to the condition $(Gg) (x) \in L_{1}((0,1);  x^{-1} dx) \cap L_{1}((1, \infty);  dx) $ and the absolute and uniform convergence with respect to $\varepsilon \in [0, 1].$  In fact, since  from (3.12) 
$$\left| \Psi_{\varepsilon}  (x,y) \right| \le  {|\Gamma (\varepsilon +ix)|^2\over \Gamma(\varepsilon + 1/2)} +
  {y^\varepsilon   \over \sqrt\pi} \int_0^\infty {dt\over (1+t)^2} \  \int_0^1   \exp \left( - {t\over u} \ -\  {yu \over 2}  \right) K_{0} \left({y u \over 2} \right)    {du\over u}  $$
$$+   {y^\varepsilon   \over \sqrt\pi} \int_0^\infty {dt\over (1+t)^2} \  \int_1^\infty    \exp \left( - {t\over u} \ -\  {yu \over 2}  \right) K_{0} \left({y u \over 2} \right)   du $$ 

$$\le  {|\Gamma (\varepsilon +ix)|^2\over \Gamma(\varepsilon + 1/2)} +   {2 y^\varepsilon   \over \sqrt\pi} \int_0^\infty { K^2_{0}\left(\sqrt{yt} \right) dt\over (1+t)^2}$$
$$  +   {y^{\varepsilon -1}  \over  \sqrt\pi} \int_0^\infty {dt\over  (1+t)^2} \  \int_0^\infty    \exp \left(-   {u \over 2}  \right) K_{0} \left({ u \over 2} \right)  \  du $$ 
$$=  O(1)+ O\left( y^{\varepsilon -1} \right),$$
we have 
$$ \int_0^\infty \left|(Gg)(y) \Psi_{\varepsilon}  (x,y)\right|   {dy\over y} \le C_1   \int_0^1  \left|(Gg)(y)\right| {dy\over y} +  C_2 \int_1^\infty \left|(Gg)(y)\right| dy < \infty,$$ 
where $C_1, C_2 >0$ are absolute constants.  Hence, combining with (3.9) and (2.21) we arrive at the inversion formula (3.5). Theorem 4 is proved.  

 \end{proof}

\section{Initial   value problem}

In this section we will apply the index transform (1.2) to investigate  the  solvability  of an initial value  problem  for the following third   order partial differential  equation, involving the Laplacian
$$\left(  x {\partial  \over \partial x}  + y {\partial   \over \partial y}  + 2 \right)  \Delta u + {\sqrt{ x^2+y^2} - 8\over x^2+ y^2} \left[ x {\partial u \over \partial x}  + y {\partial  u \over \partial y}\right]   +  {u\over 2 \sqrt{ x^2+y^2}}   =0, \  (x,y) \in 
\mathbb{R}^2 \backslash \{0\},\eqno(4.1)$$ 
where $\Delta = {\partial^2 \over \partial x^2} +  {\partial^2 \over \partial y ^2}$ is the Laplacian in $\mathbb{R}^2$.   In fact, writing  (4.1) in polar coordinates $(r,\theta)$, we end up with the equation  
$$ r { \partial^3 u \over \partial r^3} + {1\over  r} \  { \partial^3 u \over \partial r \  \partial^2 \theta} +  3 { \partial^2 u \over \partial r^2} +   \left[ 1-  {7\over r}  \right]   { \partial u \over \partial r} +  {u \over 2r}   = 0.\eqno(4.2)$$

{\bf Lemma 3.} {\it  Let $g(\tau)  \in L_1\left(\mathbb{R}; e^{ \beta  |\tau|} d\tau\right),\  \beta \in (0, 2\pi)$. Then  the function
$$u(r,\theta)=  \sqrt\pi \   \int_{-\infty}^\infty   {\rm Re}  \left[ J^2_{i\tau} \left(\sqrt  r \right)\right] \   { e^{\theta \tau} g(\tau) d\tau\over \cosh(\pi\tau)},\eqno(4.3)$$
 satisfies   the partial  differential  equation $(4.2)$ on the wedge  $(r,\theta): r   >0, \  0\le \theta <  \beta$, vanishing at infinity.}

\begin{proof} The proof  is straightforward by  substitution (4.3) into (4.2) and the use of  (1.15).  The necessary  differentiation  with respect to $r$ and $\theta$ under the integral sign is allowed via the absolute and uniform convergence, which can be verified  using  inequality (3.3)  and the integrability condition $g \in L_1\left(\mathbb{R}; e^{ \beta |\tau|} d\tau\right),\  \beta \in (0, 2\pi)$ of the lemma.  Finally,  the condition $ u(r,\theta) \to 0,\ r \to \infty$  is due to the asymptotic formula (1.4) for Bessel function at infinity. 
\end{proof}

Finally  we will formulate the initial value problem for equation (4.2) and give its solution.

{\bf Theorem 5.} {\it Let  $g(x)$ be given by formula $(3.5)$ and its transform $(Gg) (t)\equiv G(t)$ satisfies conditions of Theorem 4.  Then  $u (r,\theta),\   r >0,  \  0\le \theta < \beta$ by formula $(4.3)$  will be a solution  of the initial value problem for the partial differential  equation $(4.2)$ subject to the initial condition}
$$u(r,0) = G (r).$$

\bigskip
\centerline{{\bf Acknowledgments}}
\bigskip

The work was partially supported by CMUP (UID/MAT/00144/2013), which is funded by FCT(Portugal) with national (MEC) and European structural funds through the programs FEDER, under the partnership agreement PT2020.

\bibliographystyle{amsplain}

\begin{thebibliography}{10}

\bibitem{yak}   Yakubovich S.  Index transforms.   Singapore:  World Scientific Publishing Company; 1996.

\bibitem{erd}    Erd\'elyi A,  Magnus W,   Oberhettinger  F,   Tricomi FG.  Higher transcendental functions. Vols. I,  II. New  York: McGraw-Hill;  1953.

\bibitem{yal}   Yakubovich S,  Luchko Yu.  The hypergeometric approach to integral transforms and convolutions, Mathematics and its applications.  Vol. 287.  Dordrecht:  Kluwer Academic Publishers Group; 1994.

\bibitem{prud}  Prudnikov AP,  Brychkov  YuA,  Marichev OI. Integrals and series:  Vol. I: Elementary functions. New York:  Gordon and Breach;   1986;   Vol. II:  Special functions. New York: Gordon and Breach;  1986;   Vol. III:  More special functions. New York:   Gordon and Breach; 1990.

\bibitem {tit}  Titchmarsh EC.   An introduction to the theory of Fourier integrals.   New York:  Chelsea; 1986.

\bibitem{square}  Lebedev  NN.   On an integral representation of an arbitrary function in terms of squares of Macdonald functions with imaginary index,  {\it Sibirsk. Mat. Zh.},   {\bf 3}  (1962),  213- 222 (in Russian).

\bibitem{skal}  Yakubovich S.  New index transforms of the Lebedev- Skalskaya type, arXiv: 1509.02764.







\end{thebibliography}

\end{document}